\documentclass[12pt]{article}
\usepackage{amsmath,amsfonts,amssymb}
\usepackage{amssymb}
\usepackage{epsfig}
\author{M.~Af\/fouf\\
Department of Mathematics\\
Kean University\\
Union, NJ 07083 \\}
\title{ENERGY ESTIMATES FOR SOLUTIONS OF SPINODAL DECOMPOSITION PROBLEM}
\newcommand{\pa}{\partial}

\begin{document}
\maketitle

\vspace{1 in}
\begin{abstract}
We show the global existence of smooth solutions of a nonlinear partial differential
equation modeling the dynamics of spinodal decomposition in diffusive
materials.
{\bf AMS(MOS) Subject classifications.} Primary 35G30, 35L65.
\end{abstract}
\section{Introduction}
In this paper, we consider an equation modeling phase separation in
 spinodal decomposition dynamics, which takes place in solid and
liquid solutions under specific thermodynamical conditions.
The initial stages of phase separation is revealed in the traditional
Cahn-Hilliard theory \cite{8}. However, Aifantis and Serrin in \cite{3} suggested
a generalization of the Cahn-Hilliard theory by including additional
terms of interfacial stress. The derivation of their model equations is based
on the balance laws of mass and momentum 
\begin{eqnarray}\label{1.1}
u_t +\triangledown  J&=&0 \nonumber \\
\triangledown T&=&F 
\end{eqnarray}
where $u$ is the concentration, $ J$ is the flux of the diffusing material, 
$ F$ is the diffusive force and $ T $ is the symmetric stress tensor
which includes the interfacial terms of a typical liquid-vapor phase
transition. Combining the components of the tensor $ T $ in the one-dimensional
case leads to the expression
\begin{equation} \label{1.2}
T=-p(u) + \varepsilon u_x^2 +\delta u_{xx}
\end{equation} 
where $\delta $ is an interfacial coefficient and $ \varepsilon $ is
a short range deformity coefficient. Generally, these coefficients are functions
of concentration, but in this work, we assume them to be constants.
The equation of state $ p(u) $ is assumed to be nonconvex with a single loop.
On physical basis, the diffusive force $ F $ can be taken to be
proportional to the flux $ J $ and its time rate, to incorporate inertia 
effects, that is
\begin{equation} \label{1.3}
F=M^{-1} J + m J_t
\end{equation} 
where $ M $ is a mobility coefficient and $ m $ is a constant measuring
the effect of inertia, which will be dropped in the current model, to obtain
the relation
\begin{equation} \label{1.4}
J= -M( p(u) - \varepsilon u_x ^2 - \delta u_{xx} )_x
\end{equation} 
Furthermore, we assume that the time dependent tensor $ T $
contains viscous relaxation terms of the form $ \nu u_t$. 
Combining the relaxation terms and letting $ M=1$ in (\ref{1.4}) from the
one-dimensional mass balance equation (\ref{1.1})
we arrive at the equation
\begin{equation} \label{1.5}
u_t=p(u)_{xx}+\nu u_{xxt} -\varepsilon  (u_x^2)_{xx} -\delta u_{xxxx}
\end{equation}
This is a fourth order nonlinear differential equation, solutions
of which may explain
the later stages of the transient spinodal decomposition process.
Many aspects of the Cahn-Hilliard equation have been studied
by Bates and Fife \cite{4}, Temam \cite{11} and Witelski \cite{12}. The stationary
and mechanical solutions of (\ref{1.5}) have been investigated
by Aifantis and Serrin in \cite{2,3}.
In the current paper, we discuss the solvability of the equation (\ref{1.5}) for
various boundary conditions, which are needed to explain the long time behavior
and the evolution of spinodal decomposition process.

\section{Energy Estimates}
We consider the nonlinear evolution equation of spinodal decomposition for the density $u=u(x,t)$, that is
\begin{equation} \label{2.1}
u_t=p(u)_{xx}+\nu u_{xxt} -\varepsilon (u_x^2)_{xx} -\delta u_{xxxx} 
\end{equation}
on a bounded domain $\Omega =[0,l]$ with the initial condition
\begin{equation} \label{2.2}
u(0,x)=u_0(x) \in H^2(\Omega)
\end{equation}
The equation (\ref{2.1}) is supplemented with either periodic boundary conditions
 (see Temam [9]):
\begin{equation} \label{2.3}
\frac {\pa ^i u}{\pa x^i} (0,t)=\frac {\pa ^i u}{\pa x^i} (l,t)
\quad \mbox {for } i=0,1,2,3
\end{equation}
or Neumann boundary conditions, see Bates and Fife \cite{4}:
\begin{equation} \label{2.4}
\frac {\pa  u}{\pa x} (x,t)=\frac {\pa ^3 u}{\pa x^3} (x,t)=0
\quad \mbox {for }x=0,l.
\end{equation}
In addition, we will assume that the equation of state $p(u) \in C^{(3)}$ and
grows linearly for $\vert u \vert >N$
for some large positive number $N$ and $ p'(u) $ changes its sign inside
an interval of displacement (phase separation).\\

{\bf Local Existence and Uniqueness}\\
The derivation of local existence and uniqueness can be found for general differential equations
in Henry \cite{9}, also presented briefly in Zheng \cite{10}. 
Along the outline proof of local existence in time (see Henry \cite{9}), 
we partition the differential operators into auxiliary linearized
part and the remaining terms as follows:

\begin{equation} \label{10.1}
u_t=[-(1-\nu \partial _x^2)^{-1}\partial_x^2](\delta u_{xx})+
[(1-\nu \partial _x^2)^{-1}\partial_x^2](p(u)+ \varepsilon (u_x^2)) 
\end{equation}
We apply Fourier transform to the equation and treat the linear part as heat operator
 and denote the nonlinear part
 by $f$ to obtain
 \begin{equation} \label{11.1}
\hat u= e^{-\delta \frac{\xi ^4}{1+\nu \xi^2}t} \hat f(\xi)  
\end{equation}
We apply the contraction mapping theorem for sufficiently small time and energy estimation for the nonlinear part.
The uniqueness can be shown by applying the above partition to the difference of two solutions $w=u-v$ and
concluding that $w$ must be identically zero.

However the crucial step in proving the global existence is to have  
 a priori uniform estimates of the solution for
any time $T<+\infty$ followed by the continuation argument (see  Zheng \cite{10}). 

{\bf Global Existence}\\
Throughout this paper, $\Vert .\Vert$ will denote the $L^2(\Omega)$ norm and $c>0$
will denote a generic constant that might depend on the initial data,  $\varepsilon ,\,
\delta,\, \nu , $ and possibly $T$ but independent of $t$. The arguments in integrals will 
be omitted if they are clear. We prove the following theorem:
{\bf Main Theorem }{\em  
The equation (\ref{2.1}), equipped with the initial conditions (\ref{2.2}) and the boundary
conditions (\ref{2.3}) or (\ref{2.4}),
has a unique global solution $u \in C([0,T];H(\Omega))$.}

Proof. The proof of this theorem is based on establishing  a priori uniform
estimates on the solution $u$. We group these estimates into
two lemmas. 

\newtheorem{lemma}{Lemma}[section]
\begin{lemma}\label{lem2.1}
For any $t\in [0,T] $ we have
\begin{equation}\label{2.5}
\sup _{0\le t\le T}(\Vert u\Vert^2 + \Vert u_x\Vert^2) + \int_0^T \Vert u_{xx} \Vert^2 dt
\leq c.
\end{equation}
\end{lemma}
Proof.  Multiply equation (\ref{2.1}) by $ u $ and integrate by parts
over $\Omega $ to obtain: 
\begin{equation}\label{2.6}
 \int^l_0 u u_t dx =\int^l_0  up_{xx} dx -\nu \int^l_0  uu_{txx} dx -\varepsilon 
\int^l_0  u(u^2_x)_{xx} dx -\delta \int^l_0  uu_{xxxx}dx
\end{equation}
 We evaluate each term using the boundary conditions (\ref{2.3}) or (\ref{2.4}) and the 
restrictions on the pressure  $p$,  as follows
\begin{eqnarray}\label{2.7}
\int^l_0 u u_{txx}dx&=&u u_{tx}\bigg \vert _0^l -\int^l_0  u_xu_{tx} dx=-\frac {1}{2}(\| u_x\| ^2)_t.\nonumber \\
-\int^l_0  uu_{xxxx}dx &=&-uu_{xxx} \bigg \vert _0^l+\int^l_0  u_x u_{xxx}dx =\nonumber \\
&=& u_xu_{xx}\bigg \vert _0^l -\int^l_0  u_{xx}^2 dx = -\| u_{xx}\| ^2 .\\
-\int^l_0  u(u_x^2)_{xx} dx &=&-u(u_x^2)_x\bigg \vert_0^l+\int^l_0  u_xu(_x^2)_x dx=\nonumber \\
&=&-\int^l u_x^2 u_{xx} dx \nonumber =-\frac {1}{3} u_x^3 \bigg \vert _0^l =0.\nonumber \\
\bigg \vert \int^l_0  up_{xx} dx \bigg \vert &=& \bigg \vert (up_x)\vert _0^l -\int^l_0 
p^{\prime}(u) u_x^2 dx \bigg \vert \le k\int^l_0  u_x^2dx =k\| u_x\| ^2. \nonumber 
\end{eqnarray}
We combine these estimates to obtain
\begin{equation}\label{2.8}
\frac 12 (\| u\| ^2 +\nu \| u_x\|^2)_t +\delta \| u_{xx}\|^2 
\le k\| u_x \| ^2. 
\end{equation}
Integrate the inequality (\ref{2.8}) over time interval $[0,T]$ and invoke Gronwall's 
lemma to deduce the required estimate (\ref{2.5}).

\begin{lemma}\label{lem2.2} 
There holds
\begin{equation}\label{2.9} 
\sup _{0\le t\le T} (\| u_t\|^2 +\| u_{tx} \| ^2) +\int^T_0 \| u_{txx}\| dt
\leq c. 
\end{equation}
\end{lemma}
Proof. Differentiate (\ref{2.1}) with respect to $t$ and multiply by $u_t$ to get
\begin{equation}\label{2.10}
u_tu_{tt} = u_t p_{xxt} +\nu u_t u_{xxtt} -\varepsilon  u_t (u_x^2)_{xxt}
-\delta u_t u_{xxxxt}. 
\end{equation}
Integrating (\ref{2.10}) by parts over $\Omega $ taking into account the boundary conditions
and the estimates from lemma \ref{lem2.1} yields the following relations
\begin{eqnarray}\label{2.11}
\int^l_0 u_t u_{xxtt}dx &=&(u_tu_{xtt})\vert _0^l -\int^l_0 u_{xt}u_{xtt}dx= \nonumber\\
&=&-\frac 12 \int^l_0 u_{xt}^2 dx =-\frac 12 \| u_{xt}\| ^2.\nonumber\\
-\int^l_0 u_t (u_x^2)_{xxt} dx &=&(u_t(u_x^2)_{xt})\vert_0^l+\int^l u_{xt}
(u_x^2)_{xt} dx= \\
&=&-\int^l_0 u_{xxt}(u_x^2)_tdx=\nonumber\\
&=& -2\int^l_0 u_{xxt}u_xu_{xt}dx. \nonumber
\end{eqnarray}
Apply the Cauchy-Schwartz inequality to the last integral to get
\begin{eqnarray}\label{2.12}
\bigg \vert \int^l_0 u_t(u_x^2)_{xxt}dx\bigg \vert&\le& 2 \chi\int^l_0 u_{xxt}^2 +
\frac 4{\chi}\int^l_0 u_x^2u_{xt}^2dx \nonumber\\
&\le&2\chi \| u_{xxt}\|^2+\frac 8{\chi}\| u_{xt}\|^2_{L^{\infty}}
\int^l_0 u^2dx\\
&\le&2\chi \| u_{xxt}\|^2+\frac {8c}{\chi}\| u_{xt}\|^2_{L^{\infty}}
\int^l_0 u_x^2 dx.\nonumber
\end{eqnarray}
We evaluate the last term in (\ref{2.12}) by applying the
Young inequality
\begin{eqnarray}\label{2.13}
\Vert u_{xt}\Vert_{L^\infty} &\leq& c_1 \Vert u_{xxt} 
\Vert^{1/2} \Vert u_{xt}\Vert^{1/2} + c_2\Vert u_{xt}\Vert \nonumber\\
&\leq& c_1 \chi^2 \Vert u_{xxt} \Vert + \frac{c_1}{4\chi^2} 
\Vert u_{xt} \Vert + c_2 \Vert u_{xt}\Vert.
\end{eqnarray}
Regrouping these estimates to obtain 
\begin{equation}\label{2.14}
\Vert u_{xt} \Vert_{L^\infty}^2  \leq c\chi^2  \Vert u_{xxt}
  \Vert^2  + c(\chi) \Vert u_{xt} \Vert^2. 
 \end{equation}
where $c(\chi) $ is a function of $\chi$. Similarly, 
we evaluate the remaining terms to obtain
  \begin{equation}\label{2.15}
  \bigg\vert \int^l_0 u_tp(u)_{xxt} dx \bigg\vert\le\frac k{\chi}\| u_t\|^2+
  k\chi \| u_{xxt}\|^2
\end{equation}
  and
\begin{equation}\label{2.16} 
\int^l_0 u_tu_{xxxxt}dx=\| u_{xxt}\|^2
\end{equation}
Substituting these estimates and selecting a small enough $\chi $
 we arrive at the following inequality
\begin{equation}\label{2.18}
\frac 1{2} (\| u_t\|^2+\nu \| u_{xt}\|^2)_t+\delta_0 \| u_{xxt}\|^2 \le
  c_1\| u_t\|^2+c_2\| u_{xt}\|^2, 
\end{equation}
  where $ \delta _0 $ is a positive constant.  Integrating (\ref{2.18}) over $[0,T]$ 
and invoking the Gronwall's lemma. We conclude (\ref{2.9}).

\newtheorem{remark}{Remark}
\begin{remark}
We can derive additional energy estimates for higher order derivatives of the solution of equation (\ref{2.1}) by interpolation
relations and requiring suitable degree of smoothness of initial data. 
\end{remark}


\begin{thebibliography}{99}

\bibitem{1} E. Aifantis, {\em On the mechanics of phase transformations},
Phase Transformation, Edited by E. Aifantis and J. Gittus (1986), 233--289
\bibitem{2} E. Aifantis and J. Serrin, {\em The mechanical theory of fluid
interfaces and Maxwell's rule}, J. Coll. Interf. Sci., 96 (1983), 517--529
\bibitem{3} E. Aifantis and J. Serrin, {\em Equilibrium solutions in the mechanical 
theory of fluid microstructures}, 
J. Coll. Interf. Sci., 96 (1983), 530--547
\bibitem{4} P. Bates and P. Fife, {\em The dynamics of nucleation for 
Cahn-Hilliard equation}, 
SIAM J. Appl. Math., 53,No 4 (1993), 990--1008
\bibitem{5} Blömker, Dirk; Maier-Paape, Stanislaus; Wanner, Thomas. 
{\em Spinodal decomposition for the Cahn-Hilliard-Cook equation}
 Comm. Math. Phys.  223  (2001), no. 3, 553--582.

\bibitem{6} J. Cahn, {\em On spinodal decomposition}, Acta Metal., 9 (1961), 795--801
\bibitem{7} J. Cahn, {\em Spinodal decomposition}, Trans. Metall. Soc. AIME, 242 (1968), 
166--180
\bibitem{8} J. Cahn and J. Hilliard, {\em Free energy of a nonuniform system.
I. Interfacialfee energy}, J. Chem. Phys., 28 (1958), 258--268
\bibitem{9} D. Henry, {\em Geometric Theory of Semilinear Parabolic Equations},
 Number 840, Lecture Notes in Mathematics, Springer (1981)
\bibitem{10} S. Zheng, {\em Nonlinear parabolic equations and hyperbolic-parabolic
coupled systems}, Longman (1995)
\bibitem{11} R. Temam, {\em Infinite-dimensional dynamical systems in Mechanics
and Physics},  Springer-Verlag, New York Inc., 68 (1988)
\bibitem{12} T. Witelski, {\em The structure of internal layers for unstable
nonlinear diffusion equations}, Stud. Appl. Math., 96 (1996), 277--300

\end{thebibliography}
\end{document}